\magnification=1200

\def\Q{{\bf {Q}}}

\def\N{{\bf N}} 
     
\def\R{{\bf R}}

  \def\cN{{\cal {N}}}

%%%%%%%%%%%%%%%%%%%%%%%%%%%%%%%%%%%%%%%%%%%%%%%%%%%%
% typoref.tex. V : January 18, 2000. 
% Author : Anthony PHAN
% Warning : syntaxe +- LaTeX 
% Sources :
% T. Lachand--Robert, ``La Ma\^\i trise de \TeX'',
% R\'ef\'erences crois\'ees;
% latex.ltx's sources;
% and of course the \TeX book.
%%%%%%%%%%%%%%%%%%%%%%%%%%%%%%%%%%%%%%%%%%%%%%%%%%%%%
%
\catcode`@=11
%
% style (look at the behavior of \item dans \bibitem too,
% and at one ,\  in \re@dreferenceslist)
% Feel free to change: 	\bibn@me (title like ``R\'ef\'erences'')
%			\bibliographym@rk (general style)
%
\def\bibn@me{R\'ef\'erences}
\def\bibliographym@rk{\centerline{{\sc\bibn@me}}
	\sectionmark\section{\ignorespaces}{\unskip\bibn@me}
	\bigbreak\bgroup
	\ifx\ninepoint\undefined\relax\else\ninepoint\fi}
%
% Beware of the \bgroup: it will be closed by \endthebibliography
%
% \refsp@ce is the spacing command that appens between multiple
% references.
%
\let\refsp@ce=\ 
\let\bibleftm@rk=[
\let\bibrightm@rk=]
%
% if you want more space between brackets...
%\let\refsp@ce=\thinspace
%\def\bibleftm@rk{[\thinspace}
%\def\bibrightm@rk{\thinspace]}
%
% frenchy stuff
%
\def\numero{n\raise.82ex\hbox{$\fam0\scriptscriptstyle o$}~\ignorespaces}
%
% new variables
%
\newcount\equationc@unt
\newcount\bibc@unt
\newif\ifref@changes\ref@changesfalse
\newif\ifpageref@changes\ref@changesfalse
\newif\ifbib@changes\bib@changesfalse
\newif\ifref@undefined\ref@undefinedfalse
\newif\ifpageref@undefined\ref@undefinedfalse
\newif\ifbib@undefined\bib@undefinedfalse
\newwrite\@auxout
%
% mark an equation
%
\def\eqnum{\global\advance\equationc@unt by 1%
\edef\lastref{\number\equationc@unt}%
\eqno{(\lastref)}}
%
% One can reference anything, just copy the former macro
% and use it so: \machin \label{truc}
% In machin you would have defined \lastref by some number
% or any text.
%
% References macros
%
% The next macros are the core of \ref and \cite commands.
% Its first argument may be ref, pageref or bib.
%
% It is too tricky to be explained.
% (It is a bit recursive.)
% It allows using \cite or \ref or ...
% with arbitrary many arguments,
% for instance:
% \cite{knuth1,knuth2,ma pomme}
%
% First argument is always ref, pageref or bib.
%
\def\re@dreferences#1#2{{%
	\re@dreferenceslist{#1}#2,\undefined\@@}}
\def\re@dreferenceslist#1#2,#3\@@{\def\next{#2}%
	\expandafter\ifx\csname#1@@\meaning\next\endcsname\relax
	??\immediate\write16
	{Warning, #1-reference "\next" on page \the\pageno\space
	is undefined.}%
	\global\csname#1@undefinedtrue\endcsname
	\else\csname#1@@\meaning\next\endcsname\fi
	\ifx#3\undefined\relax
	\else,\refsp@ce\re@dreferenceslist{#1}#3\@@\fi}
%
% notice that the former ``,\refsp@ce'' will separate
% multiple arguments. But beware of spaces
% while defining a reference or calling for it!
%
% tricky thing: \newlabel has two arguments
% {labelname}{{\lastref}{\pageref}}
% The second argument is read as two arguments
% by \newl@bel. This was necessary to get
% a jobname.aux containing the same syntax
% LaTeX would produce and use.
%
\def\newlabel#1#2{{\def\next{#1}\newl@bel#2}}
\def\newl@bel#1#2{%
	\expandafter\xdef\csname ref@@\meaning\next\endcsname{#1}%
	\expandafter\xdef\csname pageref@@\meaning\next\endcsname{#2}}
\def\label#1{{%
	\toks0={#1}\message{ref(\lastref) \the\toks0,}%
	\ignorespaces\immediate\write\@auxout%
	{\noexpand\newlabel{\the\toks0}{{\lastref}{\the\pageno}}}%
	\def\next{#1}%
	\expandafter\ifx\csname ref@@\meaning\next\endcsname\lastref%
	\else\global\ref@changestrue\fi%
	\newlabel{#1}{{\lastref}{\the\pageno}}}}
\def\ref#1{\re@dreferences{ref}{#1}}
\def\pageref#1{\re@dreferences{pageref}{#1}}
%
% bibliography macros
%
\def\bibcite#1#2{{\def\next{#1}%
	\expandafter\xdef\csname bib@@\meaning\next\endcsname{#2}}}
\def\cite#1{\bibleftm@rk\re@dreferences{bib}{#1}\bibrightm@rk}
%
% The argument of \beginthebibliography
% is any sequence of numerals which will represent
% the maximum \item's length. If you have less than 9
% \bibitem's, this argument may be {any numeral}.
% if you have between 100 and 999 \bibitem's
% this argument may be {any three numerals},
% and so on.
%
\def\beginthebibliography#1{\bibliographym@rk
	\setbox0\hbox{\bibleftm@rk#1\bibrightm@rk\enspace}
	\parindent=\wd0
	\global\bibc@unt=0
	\def\bibitem##1{\global\advance\bibc@unt by 1
		\edef\lastref{\number\bibc@unt}
		{\toks0={##1}
		\message{bib[\lastref] \the\toks0,}%
		\immediate\write\@auxout
		{\noexpand\bibcite{\the\toks0}{\lastref}}}
		\def\next{##1}%
		\expandafter\ifx
		\csname bib@@\meaning\next\endcsname\lastref
		\else\global\bib@changestrue\fi%
		\bibcite{##1}{\lastref}
		\medbreak
		\item{\hfill\bibleftm@rk\lastref\bibrightm@rk}%
		}
	}
\def\endthebibliography{\egroup\par}
%
% THE NEXT MACRO MUST BE INCLUDED
% IN THE \BYE COMMAND. FOR INSTANCE:
%
% \catcode`@=11
% \outer\def\bye{\@closeaux
% 	\par\vfill\supereject\end}
% \catcode`@=12
%
\def\@closeaux{\closeout\@auxout
	\ifref@changes\immediate\write16
	{Warning, changes in references.}\fi
	\ifpageref@changes\immediate\write16
	{Warning, changes in page references.}\fi
	\ifbib@changes\immediate\write16
	{Warning, changes in bibliography.}\fi
	\ifref@undefined\immediate\write16
	{Warning, references undefined.}\fi
	\ifpageref@undefined\immediate\write16
	{Warning, page references undefined.}\fi
	\ifbib@undefined\immediate\write16
	{Warning, citations undefined.}\fi}
%
% initialization of jobname.aux
%
\immediate\openin\@auxout=\jobname.aux
\ifeof\@auxout \immediate\write16
  {Creating file \jobname.aux}
\immediate\closein\@auxout
\immediate\openout\@auxout=\jobname.aux
\immediate\write\@auxout {\relax}%
\immediate\closeout\@auxout
\else\immediate\closein\@auxout\fi
%
% Let's read this file and open it out
%
\input\jobname.aux
\immediate\openout\@auxout=\jobname.aux
% this file will be closed by \bye.
%
% That's all, folks!
%
\catcode`@=12
%\endinput

%
\catcode`@=11
\def\bibliographym@rk{\bgroup}
%
% \bye est modifie pour la biblio et la table des matieres
%
\outer\def\bye{ 	\par\vfill\supereject\end}

\def\house#1{\setbox1=\hbox{$\,#1\,$}%
\dimen1=\ht1 \advance\dimen1 by 2pt \dimen2=\dp1 \advance\dimen2 by 2pt
\setbox1=\hbox{\vrule height\dimen1 depth\dimen2\box1\vrule}%
\setbox1=\vbox{\hrule\box1}%
\advance\dimen1 by .4pt \ht1=\dimen1
\advance\dimen2 by .4pt \dp1=\dimen2 \box1\relax}

  \def\eps{{\varepsilon}}

\def\sm{\smallskip}  \def\noi{\noindent}

\def\build#1_#2^#3{\mathrel{\mathop{\kern 0pt#1}\limits_{#2}^{#3}}}

\def\date {le\ {\the\day}\ \ifcase\month\or janvier
\or fevrier\or mars\or avril\or mai\or juin\or juillet\or
ao\^ut\or septembre\or octobre\or novembre
\or d\'ecembre\fi\ {\oldstyle\the\year}}

\font\fivegoth=eufm5 \font\sevengoth=eufm7 \font\tengoth=eufm10

\newfam\gothfam \scriptscriptfont\gothfam=\fivegoth
\textfont\gothfam=\tengoth \scriptfont\gothfam=\sevengoth

\def\pro{\noindent {\it Proof. }}

\def\smallsquare{\vbox{\hrule\hbox{\vrule height 1 ex\kern 1 ex\vrule}\hrule}}
\def\cqfd{\hfill \smallsquare\vskip 3mm}

\def\og{\leavevmode\raise.3ex\hbox{$\scriptscriptstyle 
\langle\!\langle\,$}}
\def \fg {\leavevmode\raise.3ex\hbox{$\scriptscriptstyle 
\!\rangle\!\rangle\,\,$}}

\def\rme{{\rm e}}

\def\bfx{{\bf x}}
\def\bfX{{\bf X}}

\def\cN{{\cal N}}

\def\for{\ \ \hbox{for}\ }

\def\gammabar{{\overline \gamma}}

%%%%%%%%%%%%%%%%%%%%%%%%%%%%%%%%%%%%%%%%%%%%
\centerline{}

\vskip 8mm

\centerline{\bf On the Zeckendorf representation of smooth numbers}

\vskip 13mm

\centerline{Y{\sevenrm ANN} B{\sevenrm UGEAUD} %and Hajime K{\sevenrm{ANEKO}} 
\footnote{}{\rm
}}

{\narrower\narrower
\vskip 15mm

\proclaim Abstract. {
Among other results, we establish, in a quantitative form,  
that any sufficiently large integer  
cannot simultaneously be divisible only by very small primes and 
have very few digits in its Zeckendorf representation. 
}

}

\vskip 10mm

\centerline{\bf 1. Introduction and results}

\vskip 5mm

The following general (and left intentionally vague) question was introduced and 
discussed in \cite{Bu18,BuKa18}:

\smallskip
{\it Do there exist arbitrarily large integers 
which have only small prime factors and, at the same time, few nonzero digits in their
representation in some integer base?}

\smallskip

The expected answer is {\it no} and modest steps in this direction 
have been made in \cite{Bu18,BuKa18}, 
by using a combination of estimates for linear forms in complex and $p$-adic logarithms 
of algebraic numbers. We refer to \cite{Bu18,BuKa18} for bibliographical references. 

A similar question can be asked as well for the Zeckendorf representation \cite{Zec72} 
of integers having only small prime factors. 
Let $(F_n)_{n \ge 0}$ denote the Fibonacci sequence defined by $F_0 = 0$, 
$F_1 = 1$, and $F_{n +2} = F_{n+1} + F_n$ for $n \ge 0$. Every positive 
integer $N$ can be written uniquely as a sum
$$
N = \eps_\ell F_\ell + \eps_{\ell -1} F_{\ell - 1} + \ldots + \eps_2 F_2 + \eps_1 F_1,   
$$
with $\eps_\ell = 1$, $\eps_j$ in $\{0, 1\}$, and 
$\eps_j \eps_{j+1} = 0$ for $j =1, \ldots , \ell-1$. 
This representation of $N$ is called its Zeckendorf representation. 
The number of digits of $N$ in its Zeckendorf representation is the number of positive integers $j$ 
for which $\eps_j$ is equal to $1$. 
%Setting $k = \eps_1 + \ldots + \eps_\ell$, we say that $N$ has $k$ digits its 
Recall that $F_n = (\gamma^n - \gammabar^n)/ \sqrt{5}$ for $n \ge 0$, 
where $\gamma = (1 + \sqrt{5})/2$
and $\gammabar = (1 - \sqrt{5})/2$ is the Galois conjugate of $\gamma$. 
Since $\gamma$ is a unit, we cannot apply estimates for $p$-adic linear 
forms of logarithms as in \cite{Bu18,BuKa18}, thus the method developed in these papers cannot 
be straightforwardly adapted to this question.

The purpose of the present note is to give an alternative proof of the main results of 
\cite{BuKa18}, which rests only on estimates for linear forms in complex logarithms 
of algebraic numbers and can be 
extended to more general representations than $b$-ary representations, in particular to the
Zeckendorf representation (and some other Ostrowski representations, see Section 4
for a short discussion). 
Moreover, it allows us also to extend some results of \cite{Bu18,BuEv17}. 
Note that an argument similar to ours can be found in \cite{MaRo19}.

%Throughout this note, $b$ always denotes an integer at least equal to $2$. 
Following \cite{Bu18}, for an integer $k \ge 1$, we denote by $(F_j^{(k)})_{j \ge 1}$ 
the sequence, arranged in increasing order, of all positive integers which 
have at most $k$ digits in their Zeckendorf representation. 
Said differently, $(F_j^{(k)})_{j \ge 1}$ is the increasing sequence composed of the integers 
of the form
$$
F_{n_h} + \cdots + F_{n_1}, 
\quad n_h \ge n_{h-1} + 2 \ge n_{h-2} + 4 \ge \cdots \ge n_1 + 2 h -2, \ \ n_1 \ge 1, \ \ 1 \le h \le k. 
$$
In particular, the sequence $(F_j^{(1)})_{j \ge 1}$ is the Fibonacci sequence 
$(F_j)_{j \ge 1}$.

Let $S = \{q_1, \ldots , q_s\}$ be a finite, non-empty set of distinct prime numbers.
Let $n$ be a positive integer and write $n = A q_1^{r_1} \ldots q_s^{r_s}$, where 
$r_1, \ldots , r_s$ are non-negative integers and $A$ is an integer 
relatively prime to $q_1 \ldots q_s$. We define the $S$-part $[n]_S$ 
%and the $S$-free part 
of $n$ by 
$$
[n]_S := q_1^{r_1} \ldots q_s^{r_s}.
% \quad \hbox{and} \quad (n)_S := |M|.
$$

Our first result shows that 
there are only finitely many integers which have a given number of 
digits in their Zeckendorf representation 
and whose prime divisors belong to a given finite set.

\proclaim Theorem 1.1. 
Let $k$ be a positive integer and $\eps$ a positive real number.
Let $S$ be a finite, non-empty set of prime numbers. 
Then, we have 
$$
[F_j^{(k)}]_S <  (F_j^{(k)})^{\eps}, 
$$
for every sufficiently large integer $j$. 

The case $k=1$ of Theorem 1.1 has already been established in \cite{BuEv17}. 
Theorem 1.1 implies that, 
for any given positive integer $k$, the greatest prime factor of $F_j^{(k)}$ 
tends to infinity with $j$. However, its proof, 
based on the $p$-adic Schmidt Subspace Theorem, does not allow us 
to estimate the speed with which this greatest prime factor 
%$P[F_j^{(k)}]$  
tends to infinity with $j$. 
Fortunately, we are able to derive
such an estimate by means of the theory of linear forms in logarithms of algebraic numbers.

For a positive integer $n$, let $P[n]$ denote its greatest prime factor, 
with the convention that $P[1] = 1$. 
A positive real number $B$ being given, 
a positive integer $n$ is called $B$-smooth if $P[n] \le B$.

\proclaim Theorem 1.2. 
Let $S$ be a finite, non-empty set of prime numbers. 
Let $k \ge 1$ be an integer. 
Then, there exist effectively computable positive numbers $c_1$ and $j_1$, depending 
only on $k$ and $S$, such that
$$
[F_j^{(k)}]_S \le (F_j^{(k)})^{1 - c_1}, \quad
\hbox{for $j \ge j_1$}.
$$ 
Furthermore, for every positive real number $\eps$, there exists an effectively 
computable positive number $j_2$, depending 
only on $k$ and $\eps$, such that 
$$
P[F_j^{(k)}] > \Bigl({1\over k} - \eps\Bigr) \log \log F_j^{(k)} 
\, {\log \log \log F_j^{(k)} \over \log \log \log \log F_j^{(k)}}, \quad
\hbox{for $j > j_2$}.    
$$
In particular, there exists an effectively computable positive integer $n_0$, depending 
only on $k$ and $\eps$, such that any integer $n > n_0$ which is 
$$
\Bigl({1\over k} - \eps\Bigr) (\log \log n) {\log \log \log n \over \log \log \log \log n}\hbox{-smooth}
$$
has at least $k+1$ digits in its Zeckendorf representation.

A much stronger lower bound for $P[F_j^{(1)}]$ follows from Stewart's work \cite{Ste13}, namely
$$
P[F_j^{(1)}] > j \exp (\log j /104 \log \log j), \quad \hbox{for $j$ large enough}. 
$$

The analogues of Theorems 1.2 and 1.3 and of Corollaries 1.4 and 1.5 of \cite{BuKa18} also hold for 
the Zeckendorf representation instead of the base-$b$ representation, as a 
consequence of Lemma 3.2 below. We state as Theorems 1.3 and 1.4 the 
statements analogous to Corollaries 1.4 and 1.5 of \cite{BuKa18}.

\proclaim Theorem 1.3.  
Let $b \ge 2$ be an integer. 
There exists an effectively computable positive integer $n_0$ such that any integer $n > n_0$ 
satisfies the following three assertions. 
If $n$ is
$$
{\log \log n \over 2 \log \log \log \log n} \hbox{-smooth},
\,  \hbox{then $n$ has at least} \, \, 
\log \log \log n 
$$
digits in its Zeckendorf representation. 
If $n$ is
$$
\sqrt{ \log \log n \, {\log \log \log n \over \log \log \log \log n} \,} \hbox{-smooth}, 
\,  \hbox{then $n$ has at least} \, \, 
{1 \over 3} \, \sqrt{ \log \log n \, {\log \log \log n \over \log \log \log \log n} \, }
$$ 
digits in its Zeckendorf representation. 
If $n$ is
$$
{1 \over 2} \, \log \log \log n \, {\log \log \log \log n \over \log \log \log \log \log n} \hbox{-smooth},
\,  \hbox{then $n$ has at least} \, \, 
{\log \log n \over 2 \log \log \log n} 
$$
digits in its Zeckendorf representation.

Let $S$ be a finite, non-empty set of prime numbers. 
A rational integer is an integral $S$-unit if all its prime factors belong to $S$. 
Proceeding as in \cite{BuKa18}, we can deduce from Lemma 3.2 below a 
lower bound for the number of digits in the Zeckendorf representation
of integral $S$-units.

\proclaim Theorem 1.4. 
Let $S$ be a finite set of prime numbers. 
Then, for any positive real number $\eps$, 
there exists an effectively computable positive integer $n_0$, 
depending only on $S$ and $\eps$, such that any integral $S$-unit $n$ greater than $n_0$ 
has more than 
$$
(1-\eps){\log\log n \over\log\log\log n} 
$$
digits in its Zeckendorf representation.

Our method allows us to extend Theorem 1.2 of \cite{Bu18} as follows.
For a given integer $k \ge 2$, we denote by $(u_j^{(k)})_{j \ge 1}$ 
%\bfu^{(k)} := 
the sequence, arranged in increasing order, of all positive integers which are 
not divisible by $b$ and have at most $k$ nonzero digits in their $b$-ary representation. 
Said differently, $(u_j^{(k)})_{j \ge 1}$ is the ordered sequence composed of the integers 
$1, 2, \ldots , b-1$ and those of the form
$$
d_k b^{n_k} + \ldots + d_2 b^{n_2} + d_1, 
\quad n_k > \ldots > n_2 > 0, \quad
d_1, \ldots , d_k \in \{0, 1, \ldots , b-1\}, \quad
d_1 d_k \not= 0.
$$

\proclaim Theorem 1.5. 
Let $b \ge 2$ and $k \ge 2$ be integers. 
Let $S$ be a finite, non-empty set of prime numbers. 
Then, there exist effectively computable positive numbers $c_1$ and $j_1$, depending 
only on $b, k$, and $S$, such that
$$
[u_j^{(k)}]_S \le (u_j^{(k)})^{1 - c_1}, \quad
\hbox{for every $j \ge j_1$}.
$$ 

Theorem 1.5 was proved in \cite{Bu18} for $k = 2, 3$ only. 

We also get a (very slightly) weaker result than Theorem 1.1 of \cite{BuKa18}, namely, for any 
fixed $k \ge 3$, the lower bound 
$$
P[u_j^{(k)}] > \Bigl({1 \over k-1} - \eps \Bigr) \log \log u_j^{(k)} 
\, {\log \log \log u_j^{(k)} \over \log \log \log \log u_j^{(k)}}   \quad     \eqno (1.1)
%\hbox{for $j > c_6$}.  
$$
holds for any sufficiently large integer $j$. 
Note that $k-1$ is replaced by $k-2$ in Theorem 1.1 of \cite{BuKa18}. 

The paper is organized as follows. Theorem 1.1 is proved in Section 2 and Theorems 1.2 to 1.5
are established in Section 3. The final section contains some additional remarks.

\vskip 5mm

\goodbreak

\centerline{\bf  2. Proof of Theorem 1.1}

\vskip 5mm

Let $K$ be an algebraic number field. 
Denote by $M_K$ the set of places
of $K$. For $v$ in $M_K$, we choose a normalized absolute value $|\cdot |_v$
such that if $v$ is an infinite place, then 
$$
|x|_v=|x|^{[K_v:\R ]/[K:\Q ]}, \for x\in\Q,
$$
while if $v$ is finite and lies above the prime $p$, then
$$
|x|_v=|x|_p^{[K_v:\Q_p ]/[K:\Q ]}, \for x\in\Q. 
$$
These absolute values satisfy
the product formula 
$$
\prod_{v\in M_K} |x|_v=1, \quad \hbox{for every non-zero $x\in K$}.
$$
Moreover, if $x\in\Q$, then $\prod_{v|\infty} |x|_v=|x|$ 
and $\prod_{v|p} |x|_v=|x|_p$, where the products are taken
over all infinite places of $K$, respectively all places of $K$
lying above the prime number~$p$.

For later use, the height $h(x)$ of a non-zero $x$ in $K$ is defined by 
$$
h(x) = \sum_{v\in M_K} \, \log \max \{1, |x|_v\}.    \eqno (2.1)
$$

Let $T$ be a finite set of places of $K$, containing all the infinite places.
Define the ring of $T$-integers of $K$ by
$$
O_T:=\{ x\in K:\, |x|_v\leq 1\for v\in M_K\setminus T\}. 
$$
Further define
$$
H_T(x_1, \ldots ,  x_n):=\prod_{v\in T}\max \{|x_1|_v, \ldots ,  |x_n|_v\}, 
\for x_1, \ldots ,  x_n\in O_T.
$$
Our main tool is the following version of the 
$p$-adic Subspace Theorem, adapted from Theorem 3.1.3 of \cite{EvGy15}.

\proclaim Theorem 2.1. 
Let $K$ be an algebraic number field. 
Let $S$ be a finite set of prime numbers. 
Let $T$ be the finite set of places of $K$ composed of all the infinite places and all 
the places lying above the primes in $S$. 
For $v$ in $T$, let $L_{1v}, \ldots , L_{nv}$ be linearly independent linear forms in 
$X_1, \ldots , X_n$ with coefficients in $K$. 
Let $\eps > 0$. Then the set of solutions of 
$$
\prod_{v\in T} | L_{1v} (\bfx) \cdots  L_{nv} (\bfx) |_v
\le H_T(\bfx)^{-\eps}, \quad \hbox{in $\bfx \in O_K^n \setminus \{{\bf 0}\}$}, 
$$
is contained in a union of finitely many proper linear subspaces of $K^n$.

Now we proceed with the proof of Theorem 1.1. 

Let $k \ge 2$ be an integer and $\eps$ a positive real number.
Let $\cN$ be the set of integer $k$-tuples 
$(n_k, \ldots , n_1)$ such that 
$n_j - n_{j-1} \ge 2$ for $j = 2, \ldots , k$, $n_1 \ge 1$, and 
$$
[F_{n_k} + \cdots + F_{n_1}]_S > 
(F_{n_k} + \cdots + F_{n_1})^{\eps}.   \eqno (2.2)
$$

Assume that $\cN$ is infinite. 
Our aim is to apply Theorem 2.1 to get a contradiction. 
Let $(n_{k,i}, \ldots , n_{1,i})$, $i \ge 1$, denote an infinite subset of $\cN$ ordered  
such that $n_{k,i} > n_{k, i-1}$ for $i \ge 2$. 
For technical reasons, which will be clear later, we would need to assume that 
$$
\lim_{i \to + \infty} \, (n_{\ell,i} - n_{\ell-1,i}) = + \infty, \quad \ell = 2, \ldots , k,    
$$
a condition which has no reason to be satisfied. 
Let us explain how one can proceed to get a similar assumption. 

Observe that, for integers $m, \ell$ with $m > \ell > 0$, we have
$$
F_m = F_{\ell} F_{m - \ell + 1} + F_{\ell - 1} F_{m - \ell}.
$$
If there is an infinite subset $\cN_1$ of $\N$ such that 
$n_{k,i} - n_{k-1,i}$ tends to infinity when $i$ tends to infinity 
along $\cN_1$, then there is nothing 
more to do for the moment. Otherwise, there exist a positive integer $t$ and an 
infinite subset $\cN_2$ of $\N$ such that 
$n_{k,i} - n_{k-1,i} = t$ for $i$ in $\cN_2$. Then, instead of 
working with the $k$-tuple $(F_{n_{k,i}}, \ldots , F_{n_{1,i}})$, we work with the 
$(k-1)$-tuple $(G_{n_{k-1,i}}, F_{n_{k-2,i}}, \ldots , F_{n_{1,i}})$, where 
$$
G_{n_{k-1,i}} = F_t F_{n_{k-1,i} + 1} + F_{t - 1} F_{n_{k-1,i} }.
$$
Proceeding like this, we can assume that there are an integer $h$, with $1 \le h \le k$, and
an infinite subset $\cN_3$ of $\N$ such that 
%$h$-tuples $(n_{h,i}, \ldots , n_{1,i})$ such that
$
n_{h,i} > \ldots > n_{1,i}, \quad i \in \cN_3, 
$
$$
F_{n_{k,i}} + \cdots + F_{n_{1,i}} = G_{n_{h,i}} + \cdots + G_{n_{1,i}}, \quad i \in \cN_3,   \eqno (2.3) 
$$
$$
[G_{n_{h,i}} + \cdots + G_{n_{1,i}}]_S > 
(G_{n_{h,i}} + \cdots + G_{n_{1,i}})^{\eps}, \quad i \in \cN_3, 
$$ 
and
$$
\lim_{i \to + \infty} \, (n_{\ell,i} - n_{\ell-1,i}) = + \infty, \quad \ell = 2, \ldots , h,   
\quad i \in \cN_3,  \eqno (2.4)
$$
where, for $j = 1, \ldots , h$, 
we have
$$
G_{n_{j, i}} = a_j F_{n_{\ell(j), i} + 1} + b_j F_{n_{\ell(j), i}}, 
$$
for non-negative integers $a_j, b_j$ and $\ell(j)$ in $\{1, \ldots , k\}$. 

We are in position to apply Theorem 2.1. 

We work in the quadratic field $K := \Q(\sqrt{5})$. 
There are two complex embeddings, denoted by $| \cdot |_{\infty_1}$ and $| \cdot |_{\infty_2}$, numbered 
such that $| a + b \sqrt{5} |_{\infty_1} = | a + b \sqrt{5} |^{1/2}$ and 
$| a + b \sqrt{5} |_{\infty_2} = | a - b \sqrt{5} |^{1/2}$, for every rational numbers $a, b$. 
For $j = 1, \ldots , h-1$, we consider the linear forms in
$\bfX = (X_1, \ldots , X_{2h})$ defined by 
$$
L_{2j - 1, 1} (\bfX) := X_{2j-1} - X_{2j} / \sqrt{5}, \quad L_{2j - 1, 2} (\bfX) := X_{2j-1}, 
$$
and
$$
L_{2j, 1} (\bfX) := X_{2j}, \quad L_{2j, 2} (\bfX) :=  X_{2j}.
$$
Set also 
$$
L_{2h-1, 1} (\bfX) := X_{2h-1} - X_{2h} / \sqrt{5}, \quad L_{2h - 1, 2} (\bfX) := X_{2h-1}, 
$$
and 
$$
L_{2h, 1} (\bfX) := X_1 + X_3 + \ldots + X_{2h-1}, \quad 
L_{2h, 2} (\bfX) := X_2 + X_4 + \ldots + X_{2h}. 
$$
Let $T$ be the finite set of places of $K$ composed of all the infinite places and all 
the places lying above the primes in $S$. 
For every finite place $v$ in $T$, set 
$$
 L_{j, v} (\bfX) := X_j, \quad  j = 2, \ldots , 2h, \quad
L_{1, v} (\bfX) := X_1 + X_3 + \ldots + X_{2h-1}.
$$
Recall that $\gamma = (1 + \sqrt{5})/2$ and consider the points 
$$
\bfx_i = (G_{n_{1,i}}, a_1 \gamma^{n_{\ell(1),i} + 1} + b_1 \gamma^{n_{\ell(1),i}}, \ldots , 
G_{n_{h,i}}, a_h \gamma^{n_{\ell(h),i} + 1} + b_h \gamma^{n_{\ell(h),i}}), \quad i \in \cN_3. 
$$
By (2.2) and (2.3), we have 
$$
\prod_{v \in T} \prod_{j=1}^{2h} \, |L_{j,v} (\bfx_i)|_v \cdot 
\prod_{j=1}^{2h} \, |L_{j,1} (\bfx_i) L_{j,2} (\bfx_i)|_{\infty_1} 
\, |L_{j,1} (\bfx_i) L_{j,2} (\bfx_i)|_{\infty_2} 
\le H_T (\bfx_i)^{- \eps / 2},
$$
for every $i$ large enough in $\cN_3$. 

It then follows from Theorem 2.1 that there exist $t_1, \ldots , t_{2h}$ in $K$, not all zero, and an infinite
set $\cN_4$, contained in $\cN_3$, such that (1.1)
%$$
%[D_h b^{n_h} + \cdots + D_1 b^{n_1}]_S  > b^{\eps n_h},
%$$
%and
$$
t_{2h} G_{n_{h,i}} + t_{2h-1} (a_h \gamma^{n_{\ell(h),i}+1} + b_h \gamma^{n_{\ell(h),i}}) 
+ \cdots + t_{2} G_{n_{1,i}} + t_{1} (a_1 \gamma^{n_{\ell(1),i}+1} + b_1 \gamma^{n_{\ell(1),i}}) = 0,  \eqno (2.5)
$$
for every $i$ in $\cN_4$. Then, dividing (2.5) by $G_{n_{h,i}}$ and letting $i$ tend to infinity
along $\cN_4$, we deduce from (2.4) that 
$$
t_{2h} G_{n_{h,i}} + t_{2h-1} (a_h \gamma^{n_{\ell(h),i}+1} + b_h \gamma^{n_{\ell(h),i}}) = 0    \eqno (2.6)
$$ 
for infinitely many $i$ in $\cN_4$. Taking the Galois conjugate of (2.6), we obtain that 
$t_{2h} =  t_{2h-1} = 0$. Continuing like this, we get 
$t_1 = \ldots = t_{2h} = 0$, a contradiction. 
This shows that $\cN$ cannot be infinite, thus (2.2) has only finitely many solutions. 
This completes the proof of Theorem 1.1. 

\vskip 5mm

\centerline{\bf 3. Proofs of Theorems 1.2 to 1.5}

\vskip 5mm

The key tool for the proofs of Theorems 1.2 to 1.5 is the following immediate 
consequence of a theorem of Matveev \cite{Matv00}. 
The height $h$ of an algebraic number is defined in (2.1).

\proclaim Theorem 3.1.  
Let $n \ge 2$ be an integer. 
Let $\alpha_1, \ldots, \alpha_n$ be non-zero algebraic real numbers. 
Let $D$ be the degree over $\Q$ 
of a number field containing $\alpha_1, \ldots, \alpha_n$. 
Let $A_1, \ldots, A_n$ be real numbers with
$$
\log A_j \ge \max
\Bigl\{h(\alpha_j), {|\log \alpha_j| \over D}, {0.16 \over D} \Bigr\},
\qquad  1\le j \le n.
$$
Let $b_1, \ldots, b_n$ be integers and set
$$
B' = \max\Bigl\{1, \max\Bigl\{ |b_j| \ {\log A_j \over \log A_n} : 1 \le j \le n \Bigr\} \Bigr\}.
$$
Then, we have
$$
\log |\alpha_1^{b_1} \ldots \alpha_n^{b_n} - 1|  
> - 2 \times 30^{n+3} \, n^{4.5} \, D^{n+2} \, \log (\rme D)\, 
\log A_1 \ldots \log A_n  \,  \log (\rme B').      
$$

A key point in Theorem 3.1 is the presence of the factor $\log A_n$ in the denominator in the
definition of $B'$. It is crucial for getting a power saving in Theorems 1.2 and 1.5.

\vskip 5mm

We first establish Theorem 1.5 and (1.1). 

Let $b \ge 2$ be an integer. 
Below, the constants $c_1, c_2, \ldots$ are effectively computable and depend at most 
on $b$ and the constants $C_1, C_2, \ldots$ are absolute and effectively computable. 
Let $N$ be a positive integer greater than $b$ 
and $k$ the number of nonzero digits in its representation 
in base $b$. We assume that $b$ does not divide $N$, thus $k \ge 2$ and we write 
$$
N =: d_k b^{n_k} + \cdots + d_2 b^{n_2} + d_1 b^{n_1}, 
$$
where
$$
n_k > \cdots > n_2 > n_1=0, \quad
d_1, \ldots , d_k \in \{1, \ldots , b-1\}.
$$
Let $q_1, \ldots , q_s$ denote distinct prime numbers written in increasing order. 
There exist non-negative integers $r_1,\ldots,r_s$ and a positive integer $A$, 
coprime with $q_1 \cdots q_s$, such that
$$
N = A q_1^{r_1}\cdots q_s^{r_s}.
$$

In the case $A = 1$, the following lemma is similar to Lemma 3.1 of \cite{BuKa18}. 

\proclaim Lemma 3.2.  
Under the above notation, we have
$$
n_k \leq \Bigl(c_1 C_1^{s} k 
\Bigl( \prod_{i=1}^s \log q_i \Bigr) \, \log (k \log q_s)\Bigr)^{k-1} \max\{1, \log A\}.
%\eqno (3.3)
$$

\pro  
Since 
$$
\eqalign{
\Lambda_k :=
\Bigl| \Bigl(A \prod_{i=1}^s q_i^{r_i}\Bigr)d_k^{-1} b^{-n_k} -1\Bigr| 
= \Bigl| \Bigl( \prod_{i=1}^s q_i^{r_i}\Bigr)d_k^{-1} {A \over b^{n_k}}  -1\Bigr| & =
d_k^{-1} b^{-n_k}\sum_{h=1}^{k-1} d_h b^{n_h} \cr
&\leq b^{1+n_{k-1}-n_k}, \cr
}
$$
we get 
$$
\log \Lambda_k  \le - (n_k - n_{k-1}- 1) \, \log b. \eqno (3.1) 
$$
Set $A^* = \max\{A, \rme \}$. Obviously, $\Lambda_k$ is non-zero. 
Since $r_j \log q_j\le (n_k+1)\log b$ for $j=1,\ldots,s$, we deduce from Theorem 3.1 that 
$$
\log \Lambda_k \geq -c_2 C_2^s (\log q_1)\cdots (\log q_s) (\log A^*) \, \log {n_k  \over \log A^*}, 
%\eqno (3.5)
$$
thus, by (3.1), 
$$
n_k - n_{k-1} + \log A^* \le c_3 C_3^s (\log q_1)\cdots (\log q_s) (\log A^*) \, \log {n_k  \over \log A^*}.
\eqno (3.2) 
$$
Likewise, for $j = 2, \ldots , k-1$, we have
$$
\eqalign{
\Lambda_j := & 
\Bigl| \Bigl(A \prod_{i=1}^s q_i^{r_i}\Bigr)  b^{-n_j} 
(d_k b^{n_k - n_j} + \ldots + d_j)^{-1} -1\Bigr|  \cr
= & \Bigl| \Bigl( \prod_{i=1}^s q_i^{r_i}\Bigr)  b^{-n_j} 
{A \over d_k b^{n_k - n_j} + \ldots + d_j } -1\Bigr| 
 = {\sum_{h=1}^{j-1} d_h b^{n_h}   \over \sum_{h=j}^{k} d_h b^{n_h} }  
 \leq b^{1+n_{j-1}-n_k}, \cr
}
$$
thus,
$$
\log \Lambda_j  \le - (n_k - n_{j-1}- 1) \, \log b. \eqno (3.3) 
$$
Since $\Lambda_j$ is non-zero, we deduce from Theorem 3.1 that 
$$
\log \Lambda_j \geq -c_4 C_4^s (\log q_1) \cdots (\log q_s) (n_k - n_j + \log A^*) 
\, \log {n_k  \over n_k - n_j + \log A^*}. 
\eqno (3.4)
$$
Combining (3.3) and (3.4), we obtain 
$$
n_k - n_{j-1}  \leq c_5 C_5^s \Bigl( \prod_{i=1}^s \log q_i \Bigr) \, (n_k - n_j + \log A^*) 
\, \log {n_k  \over n_k - n_j + \log A}.     %\eqno (3.9)
$$
Consequently, we get 
$$
n_k - n_1  \leq n_k - n_1 + \log A^* \leq 
c_5 C_5^s \Bigl( \prod_{i=1}^s \log q_i \Bigr) \, (n_k - n_2 + \log A^*) 
\, \log {n_k  \over   \log A^*}    \eqno (3.5)
$$
and, for $j=3, \ldots , k-1$, 
$$
n_k - n_{j-1} + \log A^* \leq c_6 C_6^s \Bigl( \prod_{i=1}^s \log q_i \Bigr) \, (n_k - n_j + \log A^*) 
\, \log {n_k  \over  \log A^*}.     \eqno (3.6)
$$
The combination of (3.2), (3.5), and (3.6) then gives
$$
n_k - n_1 \leq \Bigl( c_7 C_7^s \Bigl( \prod_{i=1}^s \log q_i \Bigr) \Bigr)^{k-1} \, (\log A^*)
\Bigl(\log {n_k \over \log A^*} \Bigr)^{k-1}. 
$$
Since $n_1 = 1$, we get 
$$
\Bigl( {n_k \over \log A^*} \Bigr)^{1/(k-1)} \leq 
c_8 C_8^s \Bigl( \prod_{i=1}^s \log q_i \Bigr) \, (k-1) \, 
\Bigl(\log \Bigl( {n_k \over \log A^*} \Bigr)^{1/(k-1)} \Bigr),
$$
hence,
$$
n_k \leq \Bigl(c_9 C_9^{s} k
\Bigl( \prod_{i=1}^s \log q_i \Bigr) \, \log (k \log q_s)\Bigr)^{k-1} \, (\log A^*).
$$
%This allows us to considerably extend my earlier results on $S$-parts! 
This establishes Lemma 3.2. 
\cqfd

Since $N \le (n_k + 1) \log b$, Theorem 1.5 is a straightforward consequence of Lemma 3.2.

\medskip

Now, we consider the Zeckendorf representation of $N$
%(We can also deal with other Ostrowski representations, actually). 
and write 
$$
F_{m_k} + \ldots + F_{m_1} = N = A q_1^{r_1}\cdots q_s^{r_s},
$$
with $m_k \ge m_{k-1} + 2 \ge m_{k-2} + 4 \ge \ldots \ge m_1 + 2k - 2$ and $m_1 \ge 1$. 
Unlike for $b$-ary representations, we cannot assume that $m_1 = 1$. 

Below, the constants $c_{10}, c_{11}, \ldots$ and 
$C_{10}, C_{11}, \ldots$ are absolute and effectively computable. 

We establish the following analogue of Lemma 3.2.

\proclaim Lemma 3.3.  
Under the above notation, we have
$$
m_k \leq \Bigl(c_{10} C_{10}^{s} k 
\Bigl( \prod_{i=1}^s \log q_i \Bigr) \, \log (k \log q_s)\Bigr)^k  \max\{1, \log A\}. 
$$

\pro
We proceed as in the proof of Lemma 3.2, but we need to consider another linear form in logarithms to
show that $m_1$ cannot be too large. 
Observe that
$$
0 < | \sqrt{5} A q_1^{r_1}\cdots q_s^{r_s} - \gamma^{m_1} - \ldots - \gamma^{m_k}| 
\le |\gammabar|^{m_1} + \ldots + |\gammabar|^{m_k} \le \gamma,
$$
so
$$
\eqalign{
\Lambda := & |A q_1^{r_1}\cdots q_s^{r_s} \gamma^{- m_1} (1 + \gamma^{m_2 - m_1} + \ldots 
+ \gamma^{m_k - m_1})^{-1} - 1| \cr
= & \Bigl|  q_1^{r_1}\cdots q_s^{r_s} \gamma^{- m_1} {A \over 1 + \gamma^{m_2 - m_1} + \ldots 
+ \gamma^{m_k - m_1}}  - 1 \Bigr| < \gamma^{- m_k + 1}.  \cr}
$$
Observe that there is an absolute positive constant $C$ such that 
$$
h(1 + \gamma^{m_2 - m_1} + \ldots + \gamma^{m_k - m_1}) \le C (m_k- m_1). 
$$
Applying Theorem 3.1, this gives 
$$
m_k  \leq c_{11} C_{11}^s \Bigl( \prod_{i=1}^s \log q_i \Bigr) \, (m_k - m_1 + \log A^*) 
\, \log {m_k  \over   \log A^*}. 
$$
We consider the quantities analogous to the $\Lambda_j$'s occurring in the proof of Lemma 3.2
and, before applying Theorem 3.1, we need to check that there are nonzero. To this end, we may 
proceed as follows. Let $j$ be an integer with $1 \le j \le k$. Assume that 
$$
\sqrt{5} N = \gamma^{m_k} + \ldots + \gamma^{m_j}.      \eqno (3.7)
$$
Take the Galois conjugate to get 
$$
- \sqrt{5} N = \gammabar^{m_k} + \ldots + \gammabar^{m_j}.    \eqno (3.8)
$$
Subtracting (3.8) from (3.7), we get 
$$
2 N = F_{m_k} + \ldots + F_{m_j},
$$
a contradiction. 

Since the exact analogues of (3.2) and (3.5) hold in our context, 
we proceed as in the proof of Lemma 3.2 to get 
$$
\Bigl( {m_k \over \log A^*} \Bigr)^{1/k} \leq 
c_{12} C_{12}^s \Bigl( \prod_{i=1}^s \log q_i \Bigr) \, k \, 
\Bigl(\log \Bigl( {m_k \over \log A^*} \Bigr)^{1/k} \Bigr),
$$
and
$$
m_k \leq \Bigl(c_{13} C_{13}^{s} k
\Bigl( \prod_{i=1}^s \log q_i \Bigr) \, \log (k \log q_s)\Bigr)^k \, (\log A^*).
$$
This establishes Lemma 3.3. 
\cqfd

The first statement of Theorem 1.2 is a straightforward consequence of Lemma 3.1.
To obtain the other statements, as well as Theorems 1.3 and 1.4, we proceed 
exactly as in \cite{BuKa18}. We omit the details.

\vskip 5mm

\centerline{\bf 4. Further remarks}

\vskip 5mm

Let $\alpha$ be an irrational real number whose continued fraction expansion 
is given by $[a_0 ; a_1, a_2, \ldots]$. For $n \ge 0$, let $q_n$ be the denominator of 
the rational number $[a_0 ; a_1,  \ldots , a_n]$. Then
(see Theorem 3.9.1 of \cite{AlSh03} for a proof), every positive integer $N$ 
can be represented uniquely in the form 
$$
N = d_\ell q_\ell  + \ldots +  d_1 q_1 + d_0 q_0,
$$
where the integers $d_0, \ldots , d_\ell$ satisfy the following three condiitons:
\sm
1. $\quad 0 \le d_0 < a_0.$

2. $\quad 0 \le d_i \le a_{i+1}, \for i \ge 1.$

3. $\quad $ For $i \ge 1$, if $d_i = a_{i+1}$, then $d_{i-1} = 0$.

This representation of $N$ is called its representation in the Ostrowski $\alpha$-numeration 
system. When $\alpha$ is the Golden Ratio, this is precisely the Zeckendorf representation of $N$. 

For positive integers $b_1, \ldots , b_s$, 
Lenstra and Shallit \cite{LeSh93} established that, if $(p_n / q_n)_{n \ge 1}$ denotes the 
sequence of convergents to 
$$
\alpha = [0 ; \overline{b_{1}, \ldots , b_{s}}],
$$
then $q_{n+2s} = t q_{n+s} - (-1)^{s} q_n$, for $n \ge 1$, where $t = p_{s-1} + q_s$. 
%, where $t$ is the trace of the 
%product of matrices
%$$
%\pmatrix{
%b_1 & 1 \cr 1 & 0 \cr} \cdots \pmatrix{ b_s & 1 \cr 1 & 0  \cr}.
%$$
%Observe that
%$$
%t = K(b_1, \ldots , b_s) + K(b_2, \ldots , b_{s-1}) = p_{s-1} + q_s,
%$$
%where $p_m/q_m$ is the $m$-th convergent to $\beta$. 
This shows that, for $h= 0, \ldots , s-1$, the sequence $(q_{ns + h})_{n \ge 0}$ satisfies a
binary recurrence relation. More precisely, there exist $a_h$ in the quadratic field 
generated by $(t + \sqrt{t^2 - 4 \times (-1)^s})/2$ such that 
$$
q_{ns+h} = a_h \, \biggl( {t + \sqrt{t^2 - 4 \times (-1)^s} \over 2} \, \biggr)^n + 
{\overline a_h} \biggl( {t - \sqrt{t^2 - 4 \times (-1)^s} \over 2} \, \biggr)^n, \quad n \ge 0,
$$
where the bar denotes the Galois conjugate. 
A similar statement holds for every quadratic real number $\alpha$. 

More generally, Theorems 1.1 to 1.4 are still valid when the Zeckendorf representation is replaced 
by the representation in the Ostrowski $\alpha$-numeration 
system, provided that $\alpha$ is a quadratic real number.

\medskip

We have focused our attention to digital representations of integers. 
Actually, we can more generally consider finite sums of values of a given non-degenerate 
recurrence sequence of integers having a dominant root. Let
$$
u_n = f_1(n) \alpha_1^n + \ldots + f_t(n) \alpha_t^n, \quad \hbox{$n \ge 0$},   
$$
be a nondegenerate recurrence sequence of integers having a dominant root. Here, 
$f_1, \ldots , f_t$ are polynomials with coefficients in the 
algebraic number field $K := \Q(\alpha_1, \ldots , \alpha_t)$, for all integers 
$i, j$ with $1 \le i < j \le t$, the quotient $\alpha_i / \alpha_j$ is not a root of unity, and 
$|\alpha_1| > |\alpha_j|$, for $j = 2, \ldots , t$. For example, the Fibonacci sequence 
has these properties. 

For an integer $k \ge 1$, we denote by $(U_j^{(k)})_{j \ge 1}$ 
the sequence, arranged in increasing order, of all the integers 
of the form
$$
u_{n_k} + \cdots + u_{n_1}, 
\quad n_k > \ldots > n_1. 
$$
Let $S$ be a finite, non-empty set of prime numbers. 
We can proceed as in the proof of Theorem 1.2 in order to bound $[U_j^{(k)}]_S$
from above. A difficulty arises since we have to ensure that the 
quantities analogous to the $\Lambda_j$'s are nonzero. This may 
require additional assumption on the sequence $(u_n)_{n \ge 0}$. We leave the 
details to an interested reader.

\vskip 12mm

\centerline{\bf References}

\vskip 7mm

\beginthebibliography{999}

\medskip 

\bibitem{AlSh03}
J.-P. Allouche and J. Shallit, 
Automatic Sequences: Theory, Applications, Generalizations,  
Cambridge University Press, 2003.

\bibitem{Bu18}
Y. Bugeaud, 
{\it On the digital representation of integers with bounded prime factors}, 
Osaka J. Math. 55 (2018), 315--324.

\bibitem{BuEv17}
Y. Bugeaud and J.-H. Evertse, 
{\it $S$-parts of terms of integer linear recurrence sequences}, 
Mathematika 63 (2017), 840--851.

\bibitem{BuKa18}
Y. Bugeaud and H. Kaneko, 
{\it On the digital representation of smooth numbers}, 
Math. Proc. Cambridge Philos. Soc.  165 (2018), 533--540. 

\bibitem{EvGy15}
J.-H. Evertse and K. Gy\H ory,
Unit Equations in Diophantine Number Theory.
Cambridge University Press, 2015.

\bibitem{LeSh93}
H. W. Lenstra and J. O. Shallit,  
{\it Continued fractions and linear recurrences}, 
Math. Comp. 61 (1993), 351--354.

\bibitem{Matv00} 
E.\ M.\ Matveev,
{\it An explicit lower bound for a homogeneous rational linear form
in logarithms of algebraic numbers.\ II},
Izv.\ Ross.\ Acad.\ Nauk Ser.\ Mat.\  {64}  (2000),  125--180 (in Russian); 
English translation in Izv.\ Math.\  {64} (2000),  1217--1269.

\bibitem{MaRo19}
E. Mazumdar and S. S. Rout,
{\it Prime powers in sums of terms of binary recurrence sequences}, 
Monatsh. Math. 189 (2019), 695--714.

\bibitem{Ste80}
C. L. Stewart,
{\it On the representation of an integer in two different bases},
J. reine angew. Math.  319  (1980), 63--72.

\bibitem{Ste13}
C. L. Stewart,
{\it On divisors of Lucas and Lehmer numbers},
Acta Math. 211 (2013), 291--314.

\bibitem{Zec72}
E. Zeckendorf, 
{\it Repr\'esentation des nombres naturels par une somme de nombres de Fibonacci 
ou de nombres de Lucas}, 
Bull. Soc. Roy. Sci. Li\`ege 41 (1972), 179--182. 

\vskip 1cm

\noi Yann Bugeaud   \hfill{ }

\noi Institut de Recherche Math\'ematique Avanc\'ee, U.M.R. 7501

\noi Universit\'e de Strasbourg et C.N.R.S.

\noi 7, rue Ren\'e Descartes

\noi 67084 Strasbourg, FRANCE

\medskip

\noi e-mail : {\tt bugeaud@math.unistra.fr}

\bye